\documentclass[12pt,vatola]{article}

\usepackage{amsfonts}

\usepackage{mathrsfs}

\usepackage{graphicx}

\def\c{\centerline}

\def\re#1{\par\hangindent\parindent\indent\llap{#1\enspace}\ignorespaces}

\def\no{\noindent}

\begin{document}

\c{\bf\large A Generalization of}\vskip 3mm

\c{\bf\large Stokes Theorem on Combinatorial Manifolds}

\vskip 5mm

\c{Linfan Mao}\vskip 2mm

\c{\scriptsize (Chinese Academy of Mathematics and System Science,
Beijing 100080, P.R.China)}

\c{\scriptsize E-mail: maolinfan@163.com}

\vskip 5mm

\begin{minipage}{130mm}

\no{\bf Abstract}: {\small For an integer $m\geq 1$, a combinatorial
manifold $\widetilde{M}$ is defined to be a geometrical object
$\widetilde{M}$ such that for $\forall p\in\widetilde{M}$, there is
a local chart $(U_p,\varphi_p)$ enable $\varphi_p:U_p\rightarrow
B^{n_{i_1}}\bigcup B^{n_{i_2}}\bigcup\cdots\bigcup B^{n_{i_{s(p)}}}$
with $B^{n_{i_1}}\bigcap B^{n_{i_2}}\bigcap\cdots\bigcap
B^{n_{i_{s(p)}}}\not=\emptyset$, where $B^{n_{i_j}}$ is an
$n_{i_j}$-ball for integers $1\leq j\leq s(p)\leq m$. Integral
theory on these smoothly combinatorial manifolds are introduced.
Some classical results, such as those of {\it Stokes'} theorem and
{\it Gauss'} theorem are generalized to smoothly combinatorial
manifolds in this paper. }\vskip 2mm

\no{\bf Key Words}: {\small combinatorial manifold,  Stokes'
theorem, Gauss' theorem.}\vskip 2mm

\no{\bf AMS(2000)}: {\small 51M15, 53B15, 53B40, 57N16}

\end{minipage}

\vskip 6mm

\no{\bf \S $1.$ Introduction}

\vskip 4mm

\no As a localized euclidean space, an {\it $n$-manifold} $M^n$ is a
Hausdorff space $M^n$, i.e., a space that satisfies the $T_2$
separation axiom such that for $\forall p\in M^n$, there is an open
neighborhood $U_p, p\in U_p\subset M^n$ and a homeomorphism
$\varphi_p: U_p\rightarrow {\bf R}^n$. These manifolds,
particularly, differential manifolds are very important to modern
geometries and mechanics. By a notion of mathematical combinatorics,
i.e. {\it mathematics can be reconstructed from or turned into
combinatorization}([$3$]), the conception of combinatorial manifold
is introduced in $[4]$, which is a generalization of classical
manifolds and can be also endowed with a topological or differential
structure as a geometrical object.

Now for an integer $s\geq 1$, let $n_1,n_2,\cdots,n_s$ be an integer
sequence with $0< n_1< n_2<\cdots< n_s$. Choose $s$ open unit balls
$B_1^{n_1}, B_2^{n_2}, \cdots,B_s^{n_s}$, where
$\bigcap\limits_{i=1}^s B_i^{n_i}\not=\emptyset$ in ${\bf
R}^{n_1+_2+\cdots n_s}$. A {\it unit open combinatorial ball of
degree $s$} is a union

$$\widetilde{B}(n_1, n_2,\cdots,n_s) = \bigcup\limits_{i=1}^sB_i^{n_i}.$$

\no Then a combinatorial manifold $\widetilde{M}$ is defined in the
next.

\vskip 4mm

\no{\bf Definition $1.1$}\ {\it For a given integer sequence
$n_1,n_2,\cdots,n_m, m\geq 1$ with $0< n_1< n_2<\cdots< n_m$, a
combinatorial manifold $\widetilde{M}$ is a Hausdorff space such
that for any point $p\in \widetilde{M}$, there is a local chart
$(U_p,\varphi_p)$ of $p$, i.e., an open neighborhood $U_p$ of $p$ in
$\widetilde{M}$ and a homoeomorphism $\varphi_p:
U_p\rightarrow\widetilde{B}(n_1(p), n_2(p),\cdots,n_{s(p)}(p))$ with
$\{n_1(p),
n_2(p),\cdots,n_{s(p)}(p)\}\subseteq\{n_1,n_2,\cdots,n_m\}$ and
$\bigcup\limits_{p\in\widetilde{M}}\{n_1(p),
n_2(p),\cdots,n_{s(p)}(p)\}=\{n_1,n_2,\cdots,n_m\}$, denoted by
$\widetilde{M}(n_1,n_2,\cdots,n_m)$ or $\widetilde{M}$ on the
context and

$$\widetilde{{\mathcal A}}=\{(U_p,\varphi_p)|
p\in\widetilde{M}(n_1,n_2,\cdots,n_m))\}$$

\no an atlas on $\widetilde{M}(n_1,n_2,\cdots,n_m)$. The maximum
value of $s(p)$ and the dimension $\widehat{s}(p)$ of
$\bigcap\limits_{i=1}^{s(p)}B_i^{n_i}$ are called the dimension and
the intersectional dimensional of $\widetilde{M}(n_1,n_2,$
$\cdots,n_m)$ at the point $p$, respectively. }

\vskip 3mm

A combinatorial manifold $\widetilde{M}$ is called {\it finite} if
it is just combined by finite manifolds and {\it smooth} if it can
be endowed with a $C^{\infty}$ differential structure. For a
smoothly combinatorial manifold $\widetilde{M}$ and a point
$p\in\widetilde{M}$, it has been shown in $[4]$ that ${\rm
dim}T_p\widetilde{M}(n_1,n_2,\cdots,n_m) =\widehat{s}(p)+
\sum\limits_{i=1}^{s(p)}(n_i-\widehat{s}(p))$ and ${\rm
dim}T_p^*\widetilde{M}(n_1,n_2,\cdots,n_m) = \widehat{s}(p)+
\sum\limits_{i=1}^{s(p)}(n_i-\widehat{s}(p))$ with a basis

$$\{\frac{\partial}{\partial x^{hj}}|_p|1\leq
j\leq\widehat{s}(p)\}\bigcup
(\bigcup\limits_{i=1}^{s(p)}\bigcup\limits_{j=\widehat{s}(p)+1}^{n_i}
\{\frac{\partial}{\partial x^{ij}}|_p \ | \ 1\leq j\leq s\})$$

\no or

$$\{dx^{hj}|_p|\}1\leq
j\leq\widehat{s}(p)\}\bigcup
(\bigcup\limits_{i=1}^{s(p)}\bigcup\limits_{j=\widehat{s}(p)+1}^{n_i}
\{dx^{ij}|_p \ | \ 1\leq j\leq s\}$$

\no for a given integer $h,1\leq h\leq s(p)$. Denoted all $k$-forms
of $\widetilde{M}(n_1,n_2,\cdots,n_m)$ by $
\Lambda^k(\widetilde{M})$ and
$\Lambda(\widetilde{M})=\bigoplus\limits_{k=0}^{\widehat{s}(p)-s(p)\widehat{s}(p)
+\sum_{i=1}^{s(p)}n_i}\Lambda^k(\widetilde{M})$, then there is a
unique exterior differentiation
$\widetilde{d}:\Lambda(\widetilde{M})\rightarrow\Lambda(\widetilde{M})$
such that for any integer $k\geq 1$,
$\widetilde{d}(\Lambda^k)\subset\Lambda^{k+1}(\widetilde{M})$ with
conditions following hold similar to the classical tensor
analysis([1]).\vskip 3mm

($i$) $\widetilde{d}$ is linear, i.e., for $\forall \varphi,
\psi\in\Lambda(\widetilde{M})$, $\lambda\in{\bf R}$,

$$\widetilde{d}(\varphi+\lambda\psi)= \widetilde{d}\varphi\wedge\psi
+\lambda \widetilde{d}\psi$$

\no and for $\varphi\in\Lambda^k(\widetilde{M}),
\psi\in\Lambda(\widetilde{M})$,

$$\widetilde{d}(\varphi\wedge\psi)= \widetilde{d}\varphi+
(-1)^k\varphi\wedge \widetilde{d}\psi.$$

($ii$) For $f\in\Lambda^0(\widetilde{M})$, $\widetilde{d}f$ is the
differentiation of $f$.

($iii$) $\widetilde{d}^2=\widetilde{d}\cdot \widetilde{d}=0$.

($iv$) $\widetilde{d}$ is a local operator, i.e., if $U\subset
V\subset\widetilde{M}$ are open sets and $\alpha\in\Lambda^k(V)$,
then $\widetilde{d}(\alpha|_U)= (\widetilde{d}\alpha)|_U$.\vskip 2mm

Therefore, smoothly combinatorial manifolds poss a local structure
analogous smoothly manifolds. But notes that this local structure
maybe different for neighborhoods of different points. Whence,
geometries on combinatorial manifolds are {\it Smarandache}
geometries([$6$]-[$8$]).

There are two well-known theorems in classical tensor analysis,
i.e., {\it Stokes'} and {\it Gauss'} theorems for the integration of
differential $n$-forms on an $n$-manifold $M$, which enables us
knowing that

$$\int_Md\omega = \int_{\partial M}\omega$$

\no for a $\omega\in\Lambda^{n-1}(M)$ with compact supports and

$$\int_M({\rm div}X)\mu =\int_{\partial M}{\bf i}_X\mu$$

\no for a vector field $X$, where ${\bf
i}_X:\Lambda^{k+1}(M)\rightarrow\Lambda^{k}(M)$ defined by ${\bf
i}_X\varpi(X_1,X_2,\cdots,X_k)=\varpi(X,X_1,\cdots,X_k)$ for
$\varpi\in\Lambda^{k+1}(M)$. The similar local properties for
combinatorial manifolds with manifolds natural forwards the
following questions: {\it wether the Stokes' or Gauss' theorem is
still valid on smoothly combinatorial manifolds?} or if invalid,
{\it What are their modified forms for smoothly combinatorial
manifolds?}.

The main purpose of this paper is to find the revised Stokes' or
Gauss' theorem for combinatorial manifolds, namely, the Stokes' or
Gauss' theorem is still valid for $n$-forms on smoothly
combinatorial manifolds $\widetilde{M}$ if
$n\in\mathscr{H}_{\widetilde{M}}$, where
$\mathscr{H}_{\widetilde{M}}$ is an integer set determined by the
smoothly combinatorial manifold $\widetilde{M}$. For this objective,
we consider a particular case of combinatorial manifolds, i.e., the
combinatorial Euclidean spaces in the next section, then generalize
the definition of integration on manifolds to combinatorial
manifolds in Section $3$. The generalized form for Stokes' or Gauss'
theorem can be found in Section $4$. Terminologies and notations
used in this paper are standard and can be found in $[1]-[2]$ or
$[4]$ for those of manifolds and combinatorial manifolds
respectively.

\vskip 6mm

\no{\bf \S $2.$ Combinatorially Euclidean Spaces}

\vskip 4mm

\no As a simplest case of combinatorial manifolds, we characterize
combinatorially euclidean spaces of finite and generalize some
results in eucildean spaces in this section.

\vskip 4mm

\no{\bf Definition $2.1$} \ {\it For a given integer sequence
$n_1,n_2,\cdots,n_m, m\geq 1$ with $0< n_1< n_2<\cdots< n_m$, a
combinatorially eucildean space $\widetilde{\bf R}(n_1,\cdots,n_m)$
is a union of finitely euclidean spaces
$\bigcup\limits_{i=1}^{m}{\bf R}^{n_i}$ such that for $\forall
p\in\widetilde{\bf R}(n_1,\cdots,n_m)$,
$p\in\bigcap\limits_{i=1}^{m}{\bf R}^{n_i}$ with $\widehat{m}= {\rm
dim}(\bigcap\limits_{i=1}^{m}{\bf R}^{n_i})$ a constant.}

By definition, we can express a point $p$ of $\widetilde{{\bf R}}$
by an $m\times n_{m}$ coordinate matrix $[\overline{x}]$ following
with $x^{il}=\frac{x^{l}}{m}$ for $1\leq i\leq m, 1\leq
l\leq\widehat{m}$.

\[
[\overline{x}]=\left[
\begin{array}{cccccccc}
x^{11} & \cdots & x^{1\widehat{m}}
& x^{1(\widehat{m})+1)} & \cdots & x^{1n_1} & \cdots & 0 \\
x^{21} & \cdots & x^{2\widehat{m}}
& x^{2(\widehat{m}+1)} & \cdots & x^{2n_2} & \cdots & 0  \\
\cdots & \cdots & \cdots  & \cdots & \cdots & \cdots  \\
x^{m1} & \cdots & x^{m\widehat{m}} & x^{m(\widehat{m}+1)} & \cdots &
\cdots & x^{mn_{m}-1} & x^{mn_{m}}
\end{array}
\right]
\]

\vskip 3mm

For making a combinatorially Euclidean space to be a metric space,
we introduce {\it inner product of matrixes} similar to that of
vectors in the next.

\vskip 4mm

\no{\bf Definition $2.2$} \ {\it Let $(A)=(a_{ij})_{m\times n}$ and
$(B)=(b_{ij})_{m\times n}$ be two matrixes. The inner product
$\left<(A),(B)\right>$ of $(A)$ and $(B)$ is defined by}

$$\left<(A),(B)\right>=\sum\limits_{i,j}a_{ij}b_{ij}.$$

\vskip 3mm

\no{\bf Theorem $2.1$} \ {\it Let $(A),(B),(C)$ be $m\times n$
matrixes and $\alpha$ a constant. Then}

($1$) \ $\left<A,B\right>=\left<B,A\right>$;

($2$) \ $\left<A+B,C\right>=\left<A,C\right>+\left<B,C\right>$;

($3$) \ $\left<\alpha A,B\right>=\alpha\left<B,A\right>$;

($4$) \ $\left<A,A\right>\geq 0$ with equality hold if and only if
$(A)=O_{m\times n}$.

\vskip 3mm

{\it Proof} \ ($1$)-($3$) can be gotten immediately by definition.
Now calculation shows that

$$\left<A,A\right>=\sum\limits_{i,j}a_{ij}^2\geq 0$$

\no and with equality hold if and only if $a_{ij}=0$ for any
integers $i,j, 1\leq i\leq m, 1\leq j\leq n$, namely,
$(A)=O_{m\times n}$ \ \ \ $\natural$

\vskip 4mm

\no{\bf Theorem $2.2$} \ {\it $(A),(B)$ be $m\times n$ matrixes.
Then}

$$\left<(A),(B)\right>^2\leq\left<(A),(A)\right>\left<(B),(B)\right>$$

\no{\it and with equality hold only if $(A)=\lambda(B)$, where
$\lambda$ is a constant.}

\vskip 3mm

{\it Proof} \ If $(A)=\lambda(B)$, then
$\left<A,B\right>^2=\lambda^2\left<B,B\right>^2=\left<A,A\right>\left<B,B\right>$.
Now if there are no constant $\lambda$ enabling $(A)=\lambda(B)$,
then $(A)-\lambda(B)\not=O_{m\times n}$ for any real number
$\lambda$. According to Theorem $2.1$, we know that

$$\left<(A)-\lambda(B),(A)-\lambda(B)\right> >0,$$

\no i.e.,

$$\left<(A),(A)\right>-2\lambda\left<(A),(B)\right>+\lambda^2\left<(B),(B)\right> >0.$$

Therefore, we find that

$$\Delta=(-2\lambda)^2-4\left<(B),(B)\right>\geq 0,$$

\no namely,

$$\left<(A),(B)\right>^2\leq\left<(A),(A)\right>\left<(B),(B)\right>. \ \ \ \natural$$

\vskip 4mm

\no{\bf Corollary $2.1$} \ {\it For given real numbers $a_{ij},
b_{ij}$, $1\leq i\leq m,1\leq j\leq n$,}

$$(\sum\limits_{i,j}a_{ij}b_{ij})^2\leq(\sum\limits_{i,j}a_{ij}^2)
(\sum\limits_{i,j}b_{ij}^2).$$

\vskip 3mm

Let $\widetilde{O}$ be the origin of $\widetilde{\bf
R}(n_1,\cdots,n_m)$. Then $[O]=O_{m\times n_m}$. For $\forall
p,q\in\widetilde{\bf R}(n_1,\cdots,n_m)$, we also call
$\overrightarrow{Op}$ the vector correspondent to the point $p$
similar to classical euclidean space, Then
$\overrightarrow{pq}=\overrightarrow{Oq}-\overrightarrow{Op}$.
Theorem $2.2$ enables us to introduce an angle between two vectors
$\overrightarrow{pq}$ and $\overrightarrow{uv}$ for points
$p,q,u,v\in\widetilde{\bf R}(n_1,\cdots,n_m)$.

\vskip 4mm

\no{\bf Definition $2.3$} \ {\it Let $p,q,u,v\in\widetilde{\bf
R}(n_1,\cdots,n_m)$. Then the angle $\theta$ between vectors
$\overrightarrow{pq}$ and $\overrightarrow{uv}$ is determined by}

$$\cos\theta=\frac{\left<[p]-[q],[u]-[v]\right>}
{\sqrt{\left<[p]-[q],[p]-[q]\right>\left<[u]-[v],[u]-[v]\right>}}$$

\no{\it with the condition $0\leq\theta\leq\pi$.}

\vskip 3mm

\no{\bf Corollary $2.2$} \ {\it The conception of angle between two
vectors is well defined.}

\vskip 3mm

{\it Proof} \ Notice that

$$\left<[p]-[q],[u]-[v]\right>^2\leq\left<[p]-[q],[p]-[q]\right>
\left<[u]-[v],[u]-[v]\right>$$

\no by Theorem $2.2$. Thereby, we know that

$$-1\leq\frac{\left<[p]-[q],[u]-[v]\right>}
{\sqrt{\left<[p]-[q],[p]-[q]\right>\left<[u]-[v],[u]-[v]\right>}}\leq
1.$$

\no Therefore there is a unique angle $\theta$ with
$0\leq\theta\leq\pi$ enabling Definition $2.3$ hold. \ \ \
$\natural$

For two points $p,q$ in $\widetilde{\bf R}(n_1,\cdots,n_m)$, the
{\it distance} $d(p,q)$ between points $p$ and $q$ is defined to be
$\sqrt{\left<[p]-[q],[p]-[q]\right>}$. We get the following result.

\vskip 4mm

\no{\bf Theorem $2.3$} \ {\it For a given integer sequence
$n_1,n_2,\cdots,n_m, m\geq 1$ with $0< n_1< n_2<\cdots< n_m$,
$(\widetilde{\bf R}(n_1,\cdots,n_m);d)$ is a metric space.}

\vskip 3mm

{\it Proof} \ We only need to verify each condition for a metric
space is hold in $(\widetilde{\bf R}(n_1,\cdots,n_m);d)$. For two
point $p,q\in\widetilde{\bf R}(n_1,\cdots,n_m)$, by definition we
know that

$$d(p,q)=\sqrt{\left<[p]-[q],[p]-[q]\right>}\geq 0$$

\no with equality hold if and only if $[p]=[q]$, namely, $p=q$ and

$$d(p,q)=\sqrt{\left<[p]-[q],[p]-[q]\right>}=\sqrt{\left<[q]-[p],[q]-[p]\right>}=d(q,p).$$

Now let $u\in\widetilde{\bf R}(n_1,\cdots,n_m)$. Then by Theorem
$2.2$, we find that

\begin{eqnarray*}
& \ &(d(p,u)+d(u,p))^2\\
& \ & =\left<[p]-[u],[p]-[u]\right>+
2\sqrt{\left<[p]-[u],[p]-[u]\right>\left<[u]-[q],[u]-[q]\right>}\\
&+&\left<[u]-[q],[u]-[q]\right>\\
&\geq& \left<[p]-[u],[p]-[u]\right>+\left<[p]-[u],[u]-[q]\right>
+\left<[u]-[q],[u]-[q]\right>\\
&=& \left<[p]-[q],[p]-[q]\right>=d^2(p,q).
\end{eqnarray*}

\no Whence, $d(p,u)+d(u,p)\geq d(p,q)$ and $(\widetilde{\bf
R}(n_1,\cdots,n_m);d)$ is a metric space. \ \ \ $\natural$

\vskip 6mm

\no{\bf \S $3.$ Integration on combinatorial manifolds}

\vskip 4mm

\no We generalize the integration on manifolds to combinatorial
manifolds and show it is independent on the choice of local charts
and partition of unity in this section.

\vskip 4mm

\no{\bf $3.1$ Partition of unity}

\vskip 3mm

\no{\bf Definition $3.1$} \ {\it Let $\widetilde{M}$ be a smoothly
combinatorial manifold and $\omega\in\Lambda(\widetilde{M})$. A
support set Supp$\omega$ of $\omega$ is defined by}

$${\rm Supp}\omega = \overline{\{p\in\widetilde{M}; \omega(p)\not=0\}}$$

\no{\it and say $\omega$ has compact support if ${\rm Supp}\omega$
is compact in $\widetilde{M}$. A collection of subsets $\{C_i|
i\in\widetilde{I}\}$ of $\widetilde{M}$ is called locally finite if
for each $p\in\widetilde{M}$, there is a neighborhood $U_p$ of $p$
such that $U_p\cap C_i=\emptyset$ except for finitely many indices
$i$.}

\vskip 3mm

A {\it partition of unity} on a combinatorial manifold
$\widetilde{M}$ is defined in the next.

\vskip 4mm

\no{\bf Definition $3.2$} \ {\it A partition of unity on a
combinatorial manifold $\widetilde{M}$ is a collection $\{(U_i,g_i)|
i\in\widetilde{I}\}$, where}\vskip 3mm

($1$) \ $\{U_i| i\in\widetilde{I}\}$ is a locally finite open
covering of $\widetilde{M}$;

($2$) \ $g_i\in\mathscr{X}(\widetilde{M})$, $g_i(p)\geq 0$ for
$\forall p\in\widetilde{M}$ and ${\rm supp}g_i\in U_i$ for
$i\in\widetilde{I}$;

($3$) \ For $p\in\widetilde{M}$, $\sum\limits_{i}g_i(p)=1$.

\vskip 3mm

We get the next result for the partition of unity on smoothly
combinatorial manifolds.

\vskip 4mm

\no{\bf Theorem $3.1$} \ {\it Let $\widetilde{M}$ be a smoothly
combinatorial manifold. Then $\widetilde{M}$ admits partitions of
unity.}

\vskip 3mm

{\it Proof} \ For $\forall M\in V(G[\widetilde{M}])$, since
$\widetilde{M}$ is smooth we know that $M$ is a smoothly submanifold
of $\widetilde{M}$. As a byproduct, there is a partition of unity
$\{(U_M^{\alpha},g_M^{\alpha}) | \alpha\in I_M\}$ on $M$ with
conditions following hold.

($1$) \ $\{U_M^{\alpha}| \alpha\in I_M\}$ is a locally finite open
covering of $M$;

($2$) \ $g_M^{\alpha}(p)\geq 0$ for $\forall p\in M$ and ${\rm
supp}g_M^{\alpha}\in U_M^{\alpha}$ for $\alpha\in I_M$;

($3$) \ For $p\in M$, $\sum\limits_{i}g_M^i(p)=1$.

By definition, for $\forall p\in\widetilde{M}$, there is a local
chart $(U_p,[\varphi_p])$ enable $\varphi_p:U_p\rightarrow
B^{n_{i_1}}\bigcup B^{n_{i_2}}\bigcup\cdots\bigcup B^{n_{i_{s(p)}}}$
with $B^{n_{i_1}}\bigcap B^{n_{i_2}}\bigcap\cdots\bigcap
B^{n_{i_{s(p)}}}\not=\emptyset$. Now let $U_{M_{i_1}}^{\alpha}$,
$U_{M_{i_2}}^{\alpha}$, $\cdots$, $U_{M_{i_{s(p)}}}^{\alpha}$ be
$s(p)$ open sets on manifolds $M, M\in V(G[\widetilde{M}])$ such
that

$$p\in U_p^{\alpha}=\bigcup\limits_{h=1}^{s(p)}U_{M_{i_h}}^{\alpha}. \ \ \ \ (3.1)$$

\no We define

$$\widetilde{S}(p)=\{U_p^{\alpha}| \ {\rm all \ integers} \ \alpha \ {\rm enabling \ (3.1) \ hold} \}.$$

\no Then

$$\widetilde{\mathcal {A}}=\bigcup\limits_{p\in\widetilde{M}}\widetilde{S}(p)
=\{U_p^{\alpha}| \alpha\in\widetilde{I}(p)\}$$

\no is locally finite covering of the combinatorial manifold
$\widetilde{M}$ by properties $(1)-(3)$. For $\forall
U_p^{\alpha}\in\widetilde{S}(p)$, define

$$\sigma_{U_p^{\alpha}}=\sum\limits_{s\geq 1}
\sum\limits_{\{i_1,i_2,\cdots,i_s\}\subset\{1,2,\cdots,s(p)\}}(\prod\limits_{h=1}^sg_{M_{i_h}^{\varsigma}})$$

\no and

$$g_{U_p^{\alpha}}=\frac{\sigma_{U_p^{\alpha}}}
{\sum\limits_{\widetilde{V}\in\widetilde{S}(p)}\sigma_{\widetilde{V}}}.$$

\no Then it can be checked immediately that
$\{(U_p^{\alpha},g_{U_p^{\alpha}})|p\in\widetilde{M},
\alpha\in\widetilde{I}(p)\}$ is a partition of unity on
$\widetilde{M}$ by properties ($1$)-($3$) on $g_{M}^{\alpha}$ and
the definition of $g_{U_p^{\alpha}}$. \ \ \ $\natural$

\vskip 4mm

\no{\bf Corollary $3.1$} \ {\it Let $\widetilde{M}$ be a smoothly
combinatorial manifold with an atlas $\widetilde{\mathcal
A}=\{(V_{\alpha},[\varphi_{\alpha}])| \alpha\in\widetilde{I}\}$ and
$t_{\alpha}$ be a $C^k$ tensor field, $k\geq 1$, of field type
$(r,s)$ defined on $V_{\alpha}$ for each $\alpha$, and assume that
there exists a partition of unity $\{(U_i,g_i)| i\in J\}$
subordinate to $\widetilde{\mathcal A}$, i.e., for $forall i\in J$,
there exists $\alpha(i)$ such that $U_i\subset V_{\alpha(i)}$. Then
for $\forall p\in\widetilde{M}$,}

$$t(p)=\sum\limits_{i}g_it_{\alpha(i)}$$

\no{\it is a $C^k$ tensor field of type $(r,s)$ on $\widetilde{M}$}

\vskip 3mm

{\it Proof} \ Since $\{U_i| i\in J\}$ is locally finite, the sum at
each point $p$ is a finite sum and $t(p)$ is a type $(r,s)$ for
every $p\in\widetilde{M}$. Notice that $t$ is $C^k$ since the local
form of $t$ in a local chart $(V_{\alpha(i)},[\varphi_{\alpha(i)}])$
is

$$\sum\limits_{j}g_it_{\alpha(j)},$$

\no where the summation taken over all indices $j$ such that
$V_{\alpha(i)}\bigcap V_{\alpha(j)}\not=\emptyset$. Those number $j$
is finite by the local finiteness. \  \ \  $\natural$

\vskip 3mm

\no{\bf $3.2$ Integration on combinatorial manifolds}

\vskip 2mm

\no First, we introduce integration on combinatorial Euclidean
spaces. Let $\widetilde{\bf R}(n_1,\cdots,n_m)$ be a combinatorially
euclidean space and

$$\tau:\widetilde{\bf R}(n_1,\cdots,n_m)\rightarrow\widetilde{\bf
R}(n_1,\cdots,n_m)$$

\no a $C^1$ differential mapping with

$$[\overline{y}]=[y^{\kappa\lambda}]_{m\times n_m}=[\tau^{\kappa\lambda}([x^{\mu\nu}])]_{m\times n_m}.$$

\no The {\it Jacobi matrix} of $f$ is defined by\vskip 3mm

$$
\frac{\partial[\overline{y}]}{\partial[\overline{x}]}=
[A_{(\kappa\lambda)(\mu\nu)}],
$$

\no where $A_{(\kappa\lambda)(\mu\nu)}=\frac{\partial
\tau^{\kappa\lambda}}{\partial x^{\mu\nu}}$.

Now let $\omega\in T_k^0(\widetilde{\bf R}(n_1,\cdots,n_m))$, a {\it
pull-back} $\tau^*\omega\in T_k^0(\widetilde{\bf
R}(n_1,\cdots,n_m))$ is defined by

$$\tau^*\omega(a_1,a_2,\cdots,a_k)=\omega(f(a_1),f(a_2),\cdots,f(a_k))$$

\no for $\forall a_1,a_2,\cdots,a_k\in\widetilde{R}$.

Denoted by $n=\sum\limits_{i=1}^mn_i-\widehat{m}m$. If $0\leq l\leq
n$, recall([4]) that the basis of $\Lambda^l(\widetilde{\bf
R}(n_1,\cdots,n_m))$ is

$$\{{\bf e}^{i_1}\wedge{\bf e}^{i_2}\wedge\cdots
\wedge{\bf e}^{i_l}| 1\leq i_1<i_2\cdots<i_l\leq n\}$$

\no for a basis ${\bf e}_1,{\bf e}_2,\cdots,{\bf e}_n$ of
$\widetilde{\bf R}(n_1,\cdots,n_m)$ and its dual basis ${\bf
e}^1,{\bf e}^2,\cdots,{\bf e}^n$. Thereby the dimension of
$\Lambda^l(\widetilde{\bf R}(n_1,\cdots,n_m))$ is

\[
\left(\begin{array}{c} n\\
l
\end{array}\right)=
\frac{(\sum\limits_{i=1}^mn_i-\widehat{m}m)!}{l!(\sum\limits_{i=1}^mn_i-\widehat{m}m-l)!}.
\]

\no Whence $\Lambda^n(\widetilde{\bf R}(n_1,\cdots,n_m))$ is
one-dimensional. Now if $\omega_0$ is a basis of
$\Lambda^n(\widetilde{R})$, we then know that its each element
$\omega$ can be represented by $\omega=c\omega_0$ for a number
$c\in{\bf R}$. Let $\tau:\widetilde{\bf
R}(n_1,\cdots,n_m)\rightarrow\widetilde{\bf R}(n_1,\cdots,n_m)$ be a
linear mapping. Then

$$\tau^*:\Lambda^n(\widetilde{\bf R}(n_1,\cdots,n_m))
\rightarrow\Lambda^n(\widetilde{\bf R}(n_1,\cdots,n_m))$$

\no is also a linear mapping with
$\tau^*\omega=c\tau^*\omega_0=b\omega$ for a unique constant $b={\rm
det}\tau$, called the determinant of $\tau$. It has been known that
([$1$])

$${\rm det}\tau={\rm det}(\frac{\partial[\overline{y}]}{\partial[\overline{x}]})$$

\no for a given basis ${\bf e}_1,{\bf e}_2,\cdots,{\bf e}_n$ of
$\widetilde{\bf R}(n_1,\cdots,n_m)$ and its dual basis ${\bf
e}^1,{\bf e}^2,\cdots,{\bf e}^n$, where
$n=\sum\limits_{i=1}^mn_i-\widehat{m}m$.

\vskip 4mm

\no{\bf Definition $3.3$} \ {\it Let $\widetilde{\bf R}(n_1,
n_2,\cdots,n_m)$ be a combinatorial Euclidean
space,$n=\widehat{m}+\sum\limits_{i=1}^m(n_i-\widehat{m})$,
$\widetilde{U}\subset\widetilde{\bf R}(n_1, n_2,\cdots,n_m)$ and
$\omega\in\Lambda^n(U)$ have compact support with

$$\omega(x)=\omega_{(\mu_{i_1}\nu_{i_1})\cdots
(\mu_{i_n}\nu_{i_n})}dx^{\mu_{i_1}\nu_{i_1}}\wedge\cdots\wedge
dx^{\mu_{i_n}\nu_{i_n}}$$

\no relative to the standard basis ${\bf e}^{\mu\nu}, 1\leq\mu\leq
m, 1\leq\nu\leq n_m$ of $\widetilde{\bf R}(n_1, n_2,\cdots,n_m)$
with ${\bf e}^{\mu\nu}=e^{\nu}$ for $1\leq\mu\leq\widehat{m}$. An
integral of $\omega$ on $\widetilde{U}$ is defined to be a mapping
$\int_{\widetilde{U}}:f\rightarrow \int_{\widetilde{U}}f\in{\bf R}$
with}

$$\int_{\widetilde{U}}\omega=
\int\omega(x)\prod\limits_{\nu=1}^{\widehat{m}}dx^{\nu}\prod
\limits_{\mu\geq\widehat{m}+1,1\leq \nu\leq n_i} dx^{\mu\nu}, \ \ \
\ (3.2)$$

\no{\it where the right hand side of $(3.2)$ is the Riemannian
integral of $\omega$ on $\widetilde{U}$.}

\vskip 3mm

For example, consider the combinatorial Euclidean space
$\widetilde{\bf R}(3,5)$ with ${\bf R}^3\cap{\bf R}^5={\bf R}$. Then
the integration of an $\omega\in\Lambda^6(\widetilde{U})$ for an
open subset $\widetilde{U}\in\widetilde{\bf R}(3,5)$ is

$$\int_{\widetilde{U}}\omega=\int_{\widetilde{U}\cap({\bf R}^3\cup{\bf R}^5)}
\omega(x) dx^1dx^{12}dx^{13}dx^{22}dx^{23}dx^{24}dx^{25}.$$

\vskip 4mm

\no{\bf Theorem $3.2$} \ {\it Let $U$ and $V$ be open subsets of
$\widetilde{\bf R}(n_1,\cdots,n_m)$ and $\tau: U\rightarrow V$ is an
orientation-preserving diffeomorphism. If $\omega\in\Lambda^n(V)$
has compact support for $n=\sum\limits_{i=1}^mn_i-\widehat{m}m$,
then $\tau^*\omega\in\Lambda^n(U)$ has compact support and}

$$\int\tau^*\omega = \int\omega.$$

\vskip 3mm

{\it Proof} \ Let $\omega(x)=\omega_{(\mu_{i_1}\nu_{i_1})\cdots
(\mu_{i_n}\nu_{i_n})}dx^{\mu_{i_1}\nu_{i_1}}\wedge\cdots\wedge
dx^{\mu_{i_n}\nu_{i_n}}\in\Lambda^n(V)$. Since $\tau$ is a
diffeomorphism, the support of $\tau^*\omega$ is $\tau^{-1}({\rm
supp}\omega)$, which is compact by that of ${\rm supp}\omega$
compact.

By the usual change of variables formula, since
$\tau^*\omega=(\omega\circ\tau)({\rm det}\tau) \omega_0$ by
definition, where $\omega_0=dx^1\wedge\cdots\wedge
dx^{\widehat{m}}\wedge dx^{1(\widehat{m}+1)}\wedge
dx^{1(\widehat{m}+2)}\wedge\cdots\wedge dx^{1n_1}\wedge\cdots\wedge
dx^{mn_m}$, we then get that

\begin{eqnarray*}
\int\tau^*\omega &=&\int(\omega\circ\tau)({\rm det}\tau)
\prod\limits_{\nu=1}^{\widehat{m}}dx^{\nu}\prod\limits_{\mu\geq\widehat{m}+1,1\leq
\nu\leq n_{\mu}} dx^{\mu\nu}\\
&=& \int\omega. \ \ \ \natural
\end{eqnarray*}

\vskip 3mm

\no{\bf Definition $3.4$} \ {\it Let $\widetilde{M}$ be a smoothly
combinatorial manifold. If there exists a family
$\{(U_{\alpha},[\varphi_{\alpha}]| \alpha\in\widetilde{I})\}$ of
local charts such that}

\vskip 3mm

($1$) \
$\bigcup\limits_{\alpha\in\widetilde{I}}U_{\alpha}=\widetilde{M}$;

($2$) \ {\it for $\forall \alpha,\beta\in\widetilde{I}$, either
$U_{\alpha}\bigcap U_{\beta}=\emptyset$ or $U_{\alpha}\bigcap
U_{\beta}\not=\emptyset$ but for $\forall p\in U_{\alpha}\bigcap
U_{\beta}$, the Jacobi matrix}

$$det(\frac{\partial[\varphi_{\beta}]}{\partial[\varphi_{\alpha}]})>0,$$

\no{\it then $\widetilde{M}$ is called an oriently combinatorial
manifold and $(U_{\alpha},[\varphi_{\alpha}])$ an oriented chart for
$\forall \alpha\in\widetilde{I}$.}

For a smoothly combinatorial manifold
$\widetilde{M}(n_1,\cdots,n_m)$, it must be finite by definition.
Whence, there exists an atlas
$\mathscr{C}=\{(\widetilde{U}_{\alpha},[\varphi_{\alpha}])|
\alpha\in\widetilde{I}\}$ on $\widetilde{M}(n_1,\cdots,n_m)$
consisting of positively oriented charts such that for
$\forall\alpha\in\widetilde{I}$,
$\widehat{s}(p)+\sum\limits_{i=1}^{s(p)}(n_i-\widehat{s}(p))$ is an
constant $n_{\widetilde{U}_{\alpha}}$ for $\forall p\in
\widetilde{U}_{\alpha}$. Denote such atlas on
$\widetilde{M}(n_1,\cdots,n_m)$ by $\mathscr{C}_{\widetilde{M}}$ and
an integer family
$\mathscr{H}_{\widetilde{M}}=\{n_{\widetilde{U}_{\alpha}}|
\alpha\in\widetilde{I}\}$.

Now for any integer $n\in\mathscr{H}_{\widetilde{M}}$, we can define
an integral of $n$-forms on a smoothly combinatorial manifold
$\widetilde{M}(n_1,\cdots,n_m)$.

\vskip 4mm

\no{\bf Definition $3.5$} \ {\it Let $\widetilde{M}$ be a smoothly
combinatorial manifold with orientation $\mathscr{O}$ and
$(\widetilde{U};[\varphi])$ a positively oriented chart with a
constant $n_{\widetilde{U}}$. Suppose
$\omega\in\Lambda^{n_{\widetilde{U}}}(\widetilde{M}),
\widetilde{U}\subset\widetilde{M}$ has compact support
$\widetilde{C}\subset\widetilde{U}$. Then define

$$\int_{\widetilde{C}}\omega = \int\varphi_*(\omega|_{\widetilde{U}}). \ \ \ \ (3.3)$$

Now if $\mathscr{C}_{\widetilde{M}}$ is an atlas of positively
oriented charts with an integer set $\mathscr{H}_{\widetilde{M}}$,
let
$\widetilde{P}=\{(\widetilde{U}_{\alpha},\varphi_{\alpha},g_{\alpha})
|\alpha\in\widetilde{I}\}$ be a partition of unity subordinate to
$\mathscr{C}_{\widetilde{M}}$. For $\forall
\omega\in\Lambda^n(\widetilde{M})$,
$n\in\mathscr{H}_{\widetilde{M}}$, an integral of $\omega$ on
$\widetilde{P}$ is defined by}

$$\int_{\widetilde{P}}\omega = \sum\limits_{\alpha\in\widetilde{I}}\int g_{\alpha}\omega.
 \ \ \ \ \ (3.4)$$

\vskip 3mm

The next result shows that the integral of $n$-forms,
$n\in\mathscr{H}_{\widetilde{M}}$ is well-defined.

\vskip 4mm

\no{\bf Theorem $3.3$} \ {\it Let $\widetilde{M}(n_1,\cdots,n_m)$ be
a smoothly combinatorial manifold. For
$n\in\mathscr{H}_{\widetilde{M}}$, the integral of $n$-forms on
$\widetilde{M}(n_1,\cdots,n_m)$ is well-defined, namely, the sum on
the right hand side of $(3.4)$ contains only a finite number of
nonzero terms, not dependent on the choice of
$\mathscr{C}_{\widetilde{M}}$ and if $P$ and $Q$ are two partitions
of unity subordinate to $\mathscr{C}_{\widetilde{M}}$, then}

$$\int_{\widetilde{P}}\omega = \int_{\widetilde{Q}}\omega.$$

\vskip 3mm

{\it Proof} \ By definition for any point
$p\in\widetilde{M}(n_1,\cdots,n_m)$, there is a neighborhood
$\widetilde{U}_p$ such that only a finite number of $g_{\alpha}$ are
nonzero on $\widetilde{U}_p$. Now by the compactness of ${\rm
supp}\omega$, only a finite number of such neighborhood cover ${\rm
supp}\omega$. Therefore, only a finite number of $g_{\alpha}$ are
nonzero on the union of these $\widetilde{U}_p$, namely, the sum on
the right hand side of $(3.4)$ contains only a finite number of
nonzero terms.

Notice that the integral of $n$-forms on a smoothly combinatorial
manifold $\widetilde{M}(n_1,$ $\cdots,n_m)$ is well-defined for a
local chart $\widetilde{U}$ with a constant
$n_{\widetilde{U}}=\widehat{s}(p)+\sum\limits_{i=1}^{s(p)}(n_i-\widehat{s}(p))$
for $\forall p\in\widetilde{U}\subset\widetilde{M}(n_1,\cdots,n_m)$
by $(3.3)$ and Definition $3.3$. Whence each term on the right hand
side of $(3.4)$ is well-defined. Thereby
$\int_{\widetilde{P}}\omega$ is well-defined.

Now let
$\widetilde{P}=\{(\widetilde{U}_{\alpha},\varphi_{\alpha},g_{\alpha})
|\alpha\in\widetilde{I}\}$ and
$\widetilde{Q}=\{(\widetilde{V}_{\beta},\varphi_{\beta},h_{\beta})
|\beta\in\widetilde{J}\}$ be partitions of unity subordinate to
atlas $\mathscr{C}_{\widetilde{M}}$ and
$\mathscr{C}_{\widetilde{M}}^*$ with respective integer sets
$\mathscr{H}_{\widetilde{M}}$ and $\mathscr{H}_{\widetilde{M}}^*$.
Then these functions $\{g_{\alpha}h_{\beta}\}$ satisfy
$g_{\alpha}h_{\beta}(p)=0$ except only for a finite number of index
pairs $(\alpha,\beta)$ and

$$\sum\limits_{\alpha}\sum\limits_{\beta}g_{\alpha}h_{\beta}(p)=1,
\ \ {\rm for} \ \forall p\in\widetilde{M}(n_1,\cdots,n_m).$$

Since $\sum\limits_{\beta}=1$, we then get that

\begin{eqnarray*}
\int_{\widetilde{P}}&=& \sum\limits_{\alpha}\int g_{\alpha}\omega\\
&=& \sum\limits_{\beta}\sum\limits_{\alpha}\int
h_{\beta}g_{\alpha}\omega\\
&=& \sum\limits_{\alpha}\sum\limits_{\beta}\int
g_{\alpha}h_{\beta}\omega\\
&=& \int_{\widetilde{Q}}\omega. \ \ \ \natural
\end{eqnarray*}

Now let $n_1,n_2,\cdots, n_m$ be a positive integer sequence. For
any point $p\in\widetilde{M}$, if there is a local chart
$(\widetilde{U}_p,[\varphi_p])$ such that
$[\varphi_p]:U_p\rightarrow B^{n_1}\bigcup
B^{n_2}\bigcup\cdots\bigcup B^{n_m}$ with $B^{n_1}\bigcap
B^{n_2}\bigcap\cdots\bigcap B^{n_m}\not=\emptyset$, then
$\widetilde{M}$ is called a {\it homogenously combinatorial
manifold}. Particularly, if $m=1$, a homogenously combinatorial
manifold is nothing but a manifold. We then get consequences for the
integral of
$(\widehat{m}+\sum\limits_{i=1}^m(n_i-\widehat{m}))$-forms on
$n$-manifolds.

\vskip 4mm

\no{\bf Corollary $3.2$} \ {\it The integral of
$(\widehat{m}+\sum\limits_{i=1}^m(n_i-\widehat{m}))$-forms on a
homogenously combinatorial manifold
$\widetilde{M}(n_1,n_2,\cdots,n_m)$ is well-defined, particularly,
the integral of $n$-forms on an $n$-manifold is well-defined.}

\vskip 3mm

Similar to Theorem $3.2$ for the {\it change of variables formula of
integral} in combinatorial Euclidean space, we get that of formula
in smoothly combinatorial manifolds.

\vskip 4mm

\no{\bf Theorem $3.4$} \ {\it Let $\widetilde{M}$ and
$\widetilde{N}$ be oriently combinatorial manifolds and
$\tau:\widetilde{M}\rightarrow\widetilde{N}$ an
orientation-preserving diffeomorphism. If
$\omega\in\Lambda(\widetilde{N})$ has compact support, then
$\tau^*\omega$ has compact support and}

$$\int\omega=\int\tau^*\omega.$$

\vskip 3mm

{\it Proof} \ Notice that ${\rm supp}\tau^*\omega=\tau^{-1}({\rm
supp}\omega)$. Thereby $\tau^*\omega$ has compact support since
$\omega$ has so. Now let $\{(U_i,\varphi_i)|i\in\widetilde{I}\}$ be
an atlas of positively oriented charts of $\widetilde{M}$ and
$\widetilde{P}=\{g_i|i\in\widetilde{I}\}$ a subordinate partition of
unity with constants $n_{U_i}$. Then
$\{(\tau(U_i),\varphi_i\circ\tau^{-1})|i\in\widetilde{I}\}$ is an
atlas of positively oriented charts of $\widetilde{N}$ and
$\widetilde{Q}=\{g_i\circ\tau^{-1}\}$ is a partition of unity
subordinate to the covering $\{\tau(U_i)|i\in\widetilde{I}\}$ with
constants $n_{\tau(U_i)}$. Whence, we get that

\begin{eqnarray*}
\int\tau^*\omega &=& \sum\limits_{i}\int g_i\tau^*\omega=
\sum\limits_{i}\int \varphi_{i*}(g_i\tau^*\omega)\\
&=& \sum\limits_i\int
\varphi_{i*}(\tau^{-1})_*(g_i\circ\tau^{-1})\omega\\
&=& \sum\limits_i\int
(\varphi_i\circ\tau^{-1})_*(g_i\circ\tau^{-1})\omega\\
&=& \int\omega.   \ \ \ \ \ \ \natural
\end{eqnarray*}

\vskip 6mm

\no{\bf \S $4.$ A generalization of Stokes' theorem}

\vskip 4mm

\no{\bf Definition $4.1$} \ {\it Let $\widetilde{M}$ be a smoothly
combinatorial manifold. A subset $\widetilde{D}$ of $\widetilde{M}$
is with boundary if its points can be classified into two classes
following.}\vskip 3mm

{\bf Class $1$(interior point Int$\widetilde{D}$)} {\it For $\forall
p\in{\rm Int}D$, there is a neighborhood $\widetilde{V}_p$ of $p$
enable $\widetilde{V}_p\subset\widetilde{D}$.}

{\bf Case $2$(boundary $\partial\widetilde{D}$)} {\it For $\forall
p\in\partial\widetilde{D}$, there is integers $\mu, \nu$ for a local
chart $(U_p;[\varphi_p])$ of $p$ such that $x^{\mu\nu}(p)=0$ but}

$$\widetilde{U}_p\cap \widetilde{D}=\{q| q\in U_p, x^{\kappa\lambda}\geq 0 \ for \
\forall\{\kappa,\lambda\}\not=\{\mu,\nu\}\}.$$\vskip 2mm

Then we generalize the famous {\it Stokes theorem} on manifolds in
the next.

\vskip 4mm

\no{\bf Theorem $4.1$}\ {\it Let $\widetilde{M}$ be a smoothly
combinatorial manifold with an integer set
$\mathscr{H}_{\widetilde{M}}$ and $\widetilde{D}$ a boundary subset
of $\widetilde{M}$. For $n\in\mathscr{H}_{\widetilde{M}}$ if
$\omega\in\Lambda^{n}(\widetilde{M})$ has compact support, then}

$$\int_{\widetilde{D}}d\omega = \int_{\partial \widetilde{D}}\omega$$

\no{\it with the convention $\int_{\partial \widetilde{D}}\omega=0$
while $\partial \widetilde{D}=\emptyset$.}

\vskip 3mm

{\it Proof} \ By Definition $3.5$, the integration on a smoothly
combinatorial manifold was constructed with partitions of unity
subordinate to an atlas. Let $\mathscr{C}_{\widetilde{M}}$ be an
atlas of positively oriented charts with an integer set
$\mathscr{H}_{\widetilde{M}}$ and
$\widetilde{P}=\{(\widetilde{U}_{\alpha},\varphi_{\alpha},g_{\alpha})
|\alpha\in\widetilde{I}\}$ a partition of unity subordinate to
$\mathscr{C}_{\widetilde{M}}$. Since ${\rm supp}\omega$ is compact,
we know that

$$\int_{\widetilde{D}}d\omega=\sum\limits_{\alpha\in\widetilde{I}}
\int_{\widetilde{D}}d(g_{\alpha\omega}),$$

$$\int_{\partial\widetilde{D}}\omega=\sum\limits_{\alpha\in\widetilde{I}}
\int_{\partial\widetilde{D}}g_{\alpha\omega}.$$

\no and there are only finite nonzero terms on the right hand side
of the above two formulae. Thereby, we only need to prove

$$\int_{\widetilde{D}}d(g_{\alpha}\omega) =
\int_{\partial \widetilde{D}}g_{\alpha}\omega$$

\no for $\forall\alpha\in\widetilde{I}$.

Not loss of generality we can assume that $\omega$ is an $n$-forms
on a local chart $(U,\varphi)$ with compact support. Now write

$$\omega(x)=\omega_{(\mu_{i_1}\nu_{i_1})\cdots
(\mu_{i_n}\nu_{i_n})}dx^{\mu_{i_1}\nu_{i_1}}\wedge\cdots\wedge
dx^{\mu_{i_n}\nu_{i_n}}$$

$$\omega=\sum\limits_{h=1}^n(-1)^{h-1}\omega_{\mu_{i_h}\nu_{i_h}}
dx^{\mu_{i_1}\nu_{i_1}}\wedge\cdots\wedge
\widehat{dx^{\mu_{i_h}\nu_{i_h}}}\wedge\cdots\wedge
dx^{\mu_{i_n}\nu_{i_n}},$$

\no where $\widehat{dx^{\mu_{i_h}\nu_{i_h}}}$ means that
$dx^{\mu_{i_h}\nu_{i_h}}$ is deleted, where

$$i_h\in\{1,\cdots,\widehat{n}_U, (1(\widehat{n}_U+1)),\cdots,(1n_1),
(2(\widehat{n}_U+1)),\cdots, (2n_2),\cdots,(mn_m)\}.$$

\no Then

$$d\omega=\sum\limits_{i=1}^n\frac{\partial\omega_i}{\partial x^i}dx^{\mu_{i_1}\nu_{i_1}}\wedge\cdots\wedge
dx^{\mu_{i_h}\nu_{i_h}}\wedge\cdots\wedge dx^{\mu_{i_n}\nu_{i_n}}. \
\ \ \ (4.1)$$

Consider the appearance of chart $U$. There are two cases must be
considered.

\vskip 3mm

\no{\bf Case $1$} \ \ $U\bigcap\partial U=\emptyset$

\vskip 3mm

In this case, $\int_{\partial U}\omega=0$ and $U$ is in
$\widetilde{M}\setminus\widetilde{D}$ or in {\bf
Int}$\widetilde{D}$. The former is naturally implies that
$\int_{\widetilde{D}}d(g_{\alpha}\omega)=0$. For the later, we find
that

$$\int_{\widetilde{D}}d\omega=\sum\limits_{i=1}^n\int_U
\frac{\partial\omega_i}{\partial x^i}dx^{\mu_{i_1}\nu_{i_1}}\cdots
dx^{\mu_{i_n}\nu_{i_n}}. \ \ \ \ (4.2)$$

\no Notice that $\int_{-\infty}^{+\infty}\frac{\partial\omega_i}
{\partial x^i}dx^i=0$ since $\omega_i$ has compact support. Thus
$\int_Ud\omega=0$ as desired.

\vskip 3mm

\no{\bf Case $2$} \ \ $\partial U\not=\emptyset$

\vskip 3mm

In this case we can do the same trick for each term except the last.
Without loss of generality, assume that

$$U\bigcap\widetilde{D}=\{q | q\in U, x^n(q)\geq 0\}$$

\no and

$$U\bigcap\partial\widetilde{D}=\{q | q\in U, x^n(q)=0\}.$$

\no Then we get that

\begin{eqnarray*}
\int_{\partial \widetilde{D}}\omega &=& \int_{U\cap\partial
\widetilde{D}}\omega\\
&=& \sum\limits_{h=1}^n(-1)^{h-1}\int_{U\cap\partial
\widetilde{D}}\omega_{\mu_{i_h}\nu_{i_h}}
dx^{\mu_{i_1}\nu_{i_1}}\wedge\cdots\wedge
\widehat{dx^{\mu_{i_h}\nu_{i_h}}}\wedge\cdots\wedge
dx^{\mu_{i_n}\nu_{i_n}}\\
&=& (-1)^{n-1}\int_{U\cap\partial \widetilde{D}}\omega_{\mu_n\nu_n}
dx^{\mu_{i_1}\nu_{i_1}}\wedge\cdots\wedge
dx^{\mu_{i_{n-1}}\nu_{i_{n-1}}}
\end{eqnarray*}

\no since $dx^n(q)=0$ for $q\in U\cap\partial\widetilde{D}$. Notice
that ${\bf R}^{n-1}=\partial{\bf R}_+^n$ but the usual orientation
on ${\bf R}^{n-1}$ is not the boundary orientation, whose outward
unit normal is $-{\bf e}_n=(0,\cdots,o,-1)$. Hence

$$\int_{\partial \widetilde{D}}\omega = -\int_{\partial
{\bf
R}_+^n}\omega_{\mu_n\nu_n}(x^{\mu_{i_1}\nu_{i_1}},\cdots,x^{\mu_{i_{n-1}}\nu_{i_{n-1}}},0)
dx^{\mu_{i_1}\nu_{i_1}}\cdots dx^{\mu_{i_{n-1}}\nu_{i_{n-1}}}.$$

On the other hand, by the fundamental theorem of calculus,

\begin{eqnarray*}&\ &\int_{{\bf R}^{n-1}}(\int_0^{\infty}
\frac{\partial\omega_{\mu_n\nu_n}}{\partial
x^{\mu_n\nu_n}})dx^{\mu_{i_1}\nu_{i_1}}\cdots
dx^{\mu_{i_{n-1}}\nu_{i_{n-1}}}\\
&\ & =-\int_{{\bf
R}^{n-1}}\omega_{\mu_n\nu_n}(x^{\mu_{i_1}\nu_{i_1}}\cdots
x^{\mu_{i_n}\nu_{i_n}},0)dx^{\mu_{i_1}\nu_{i_1}}\cdots
dx^{\mu_{i_{n-1}}\nu_{i_{n-1}}}.\end{eqnarray*}

\no Since $\omega_{\mu_{i_n}\nu_{i_n}}$ has compact support, thus

$$\int_U\omega=-\int_{{\bf
R}^{n-1}}\omega_{\mu_n\nu_n}(x^{\mu_{i_1}\nu_{i_1}}\cdots
x^{\mu_{i_n}\nu_{i_n}},0)dx^{\mu_{i_1}\nu_{i_1}}\cdots
dx^{\mu_{i_{n-1}}\nu_{i_{n-1}}}.$$

\no Therefore, we get that

$$\int_{\widetilde{D}}d\omega = \int_{\partial \widetilde{D}}\omega$$

\no This completes the proof. \  \ \ \ $\natural$

Corollaries following are immediately obtained by Theorem $4.1$

\vskip 4mm

\no{\bf Corollary $4.1$} \ {\it Let $\widetilde{M}$ be a smoothly
and homogenously combinatorial manifold with an integer set
$\mathscr{H}_{\widetilde{M}}$ and $\widetilde{D}$ a boundary subset
of $\widetilde{M}$. For $n\in\mathscr{H}_{\widetilde{M}}$ if
$\omega\in\Lambda^{n}(\widetilde{M})$ has compact support, then}

$$\int_{\widetilde{D}}d\omega = \int_{\partial \widetilde{D}}\omega,$$

\no{\it particularly, if $\widetilde{M}$ is nothing but a manifold,
the Stokes theorem holds.}

\vskip 3mm

\no{\bf Corollary $4.2$} \ {\it Let $\widetilde{M}$ be a smoothly
combinatorial manifold with an integer set
$\mathscr{H}_{\widetilde{M}}$. For
$n\in\mathscr{H}_{\widetilde{M}}$, if
$\omega\in\Lambda^n(\widetilde{M})$ has a compact support, then}

$$\int_{\widetilde{M}}\omega=0.$$

\vskip 3mm

Similar to the case of manifolds, we find a generalization for {\it
Gauss theorem} in the next.

\vskip 4mm

\no{\bf Theorem $4.2$} \ {\it Let $\widetilde{M}$ be a smoothly
combinatorial manifold with an integer set
$\mathscr{H}_{\widetilde{M}}$, $\widetilde{D}$ a boundary subset of
$\widetilde{M}$ and {\bf X} a vector field on $\widetilde{M}$ with
compact support. Then}

$$\int_{\widetilde{D}}({\rm div}{\bf X}){\bf v}=\int_{\partial\widetilde{D}}
{\bf i}_{\bf X}{\bf v},$$

\no{\it where ${\bf v}$ is a volume form on $\widetilde{M}$, i.e.,
nonzero elements in $\Lambda^n(\widetilde{M})$ for
$n\in\mathscr{H}_{\widetilde{M}}$.}

\vskip 3mm

{\it Proof} \ This result is also a consequence of Theorem $4.1$.
Notice that

$$({\rm div}{\bf X}){\bf v}=d{\bf i}_{\bf X}{\bf v}+{\bf i}_{\bf X}d{\bf v}
=d{\bf i}_{\bf X}.$$

\no According to Theorem $4.1$, we then get that

$$\int_{\widetilde{D}}({\rm div}{\bf X}){\bf v}=\int_{\partial\widetilde{D}}
{\bf i}_{\bf X}{\bf v}. \ \ \ \natural$$

\vskip 8mm

\no{\bf References}

\vskip 5mm

\re{[1]}R.Abraham, J.E.Marsden and T.Ratiu, {\it Manifolds, Tensor
Analysis and Application}, Addison-Wesley Publishing Company, Inc,
Reading, Mass, 1983.

\re{[2]}W.H.Chern and X.X.Li, {\it Introduction to Riemannian
Geometry}, Peking University Press, 2002.

\re{[3]}L.F.Mao, Combinatorial speculations and the combinatorial
conjecture for mathematics, in {\it Selected Papers on Mathematical
Combinatorics}(I), World Academic Union, 2006.

\re{[4]}L.F.Mao, Geometrical theory on combinatorial manifolds, {\it
arXiv: math.GM/0612760}, will also appears in {\it JP J.Geometry and
Topology}(accepted).

\re{[5]}L.F.Mao, Pseudo-Manifold Geometries with Applications,
e-print: {\it arXiv: math.} {\it GM/0610307}.

\re{[6]}L.F.Mao, {\it Automorphism Groups of Maps, Surfaces and
Smarandache Geometries}, American Research Press, 2005.

\re{[7]}L.F.Mao, {\it Smarandache multi-space theory}, Hexis,
Phoenix, AZ£¬2006.

\re{[8]}F.Smarandache, Mixed noneuclidean geometries, {\it eprint
arXiv: math/0010119}, 10/2000.

\end{document}